\documentclass[12 pt]{article}
\usepackage{amssymb,amsmath,amsfonts,amsthm}
\newcommand\BBR{{\mathbb {R}}}
\newcommand\R{{\mathbb {R}}}

\newcommand\BBN{{\mathbb {N}}}
\newcommand\N{{\mathbb {N}}}

\newtheorem {Lemma}{Lemma}[section]
\newtheorem {Theorem}{Theorem}[section]
\newtheorem {Proposition}{Proposition}[section]

\theoremstyle{definition}

\renewcommand{\P}{\ensuremath{\mathbb {P}}}
\newcommand{\E}{\mathbb{E}}

\newcommand{\p}{ { \mathbb P} }

\newcommand\beq{\begin{equation}}
\newcommand\eeq{\end{equation}}

\begin{document}

\title{Berry-Esseen type bounds for the matrix coefficients  and the spectral radius of the left random walk on $GL_d({\mathbb R})$}

\author{C. Cuny\footnote{Christophe Cuny, Univ Brest, LMBA, UMR 6205 CNRS, 6 avenue Victor Le Gorgeu, 29238 Brest}, J. Dedecker\footnote{J\'er\^ome Dedecker, Universit\'e de Paris, CNRS, MAP5, UMR 8145,
45 rue des  Saints-P\`eres,
F-75006 Paris, France.}
and
F. Merlev\`ede \footnote{Florence Merlev\`ede, LAMA,  Univ Gustave Eiffel, Univ Paris Est Cr\'eteil, UMR 8050 CNRS,  \  F-77454 Marne-La-Vall\'ee, France.}
M. Peligrad \footnote{Magda Peligrad, Department of Mathematical Sciences, University of Cincinnati, PO Box 210025, Cincinnati, Oh 45221-0025, USA.}}

\maketitle

\begin{abstract}
We give rates of convergence in the Central Limit Theorem for the coefficients and the spectral radius of the left random walk on $GL_d({\mathbb R})$, assuming the existence of an exponential or polynomial moment. 
\end{abstract}

\section{Introduction}

Let $(\varepsilon_n)_{n \geq 1}$ be independent random matrices taking values in $G= GL_d(\mathbb R)$, $d \geq 2$ (the group of invertible $d$-dimensional real matrices) with common distribution $\mu$. Let $\Vert \cdot \Vert$ be the euclidean norm on ${\mathbb R}^d$, and for every $A \in GL_d(\mathbb R)$, let $\|A\|=\sup_{x, \|x\|=1} \|A x \|$. Let also $N(g) := \max ( \Vert g \Vert , \Vert g^{-1} \Vert)$. We shall say that $\mu $ has an exponential moment if there exists $\alpha>0$ such that 
\[
\int_G (N(g) )^\alpha d \mu(g) < \infty \, , 
\]
We shall say that $\mu $ has a polynomial moment of order $p \geq 1$ if
\[
\int_G (\log N(g) )^p d \mu(g) < \infty \,  .
\]

Let $A_n= \varepsilon_n \cdots \varepsilon_1$, with the convention $A_0=$Id. It follows from Furstenberg and Kesten \cite{FK} that, if $\mu$ admits a moment of order $1$ then 
\beq \label{SL1}
\lim_{n \rightarrow \infty} \frac{1}{n}  \log \Vert A_n \Vert = \lambda_{\mu} \, \text{ ${\mathbb P}$-a.s.},
\eeq
where $ \lambda_{\mu} := \lim_{n \rightarrow \infty} n^{-1}  \E \log \Vert A_n  \Vert $ is the so-called first Lyapunov exponent.

Let now $X:= P({\mathbb R}^d)$  be the projective space of ${\mathbb R}^d $ and write ${\bar x}$ as the projection of $x \in {\mathbb R}^d -\{0\}$ to $X$.  An element $A$ of $G= GL_d(\mathbb R)$ acts on  the projective space $X$ as follows: $A \bar x = \overline{Ax}$.  Let $\Gamma_\mu$ be the closed semi-group generated by the support of $\mu$. We say that $\mu$ is proximal if $\Gamma_\mu$ contains a matrix  that admits  a unique (with multiplicity $1$) eigenvalue of maximal modulus. We say that $\mu$ is strongly irreducible if no proper union of subspaces of ${\mathbb R}^d$ is invariant by $\Gamma_\mu$. Throughout the paper, we  assume that $\mu$ is strongly irreducible and proximal.
In particular,  there exists a unique invariant measure $\nu$ on ${\mathcal B} (X)$, meaning that for any bounded measurable function $h$ from $X$ to $\mathbb R$,
\beq \label{defnu}
\int_X h(x) d \nu(x) = \int_G \int_X h( g \cdot x ) d \mu(g) d \nu(x) \, .
\eeq
Let $W_0$ be a random variable with values in the projective space $X$, independent of $(\varepsilon_n)_{n \geq 1}$  and with distribution $\nu$. By the invariance of $\nu$, we see that the sequence $(W_n:=A_n W_0)_{n \geq 1}$  is a strictly stationary Markov chain with values in $X$. Let now, for any integer $k \geq 1$, 
\beq  \label{defMCXk}
X_k := \sigma (\varepsilon_k, W_{k-1} ) - \lambda_{\mu} = \sigma (\varepsilon_k, A_{k-1} W_0 ) - \lambda_{\mu} \, , 
\eeq 
where,   for any $g \in G$ and any ${\bar x} \in X$,
$
\sigma( g , {\bar x} ) = \log (\|g x \|/\|x\|) 
$.
Note that $\sigma$ is an additive cocycle in the sense that $\sigma ( g_1 g_2,   {\bar x})  = \sigma ( g_1,    g_2 {\bar x})  +  \sigma ( g_2,  {\bar x}) $.  Consequently
\[
S_n := \sum_{k=1}^n X_k = \log \Vert  A_n V_0 \Vert  - n \lambda_{\mu}\, ,
\]
where $V_0$ is a random variable such that $\|V_0\|=1$ and $\overline V_0= W_0$.

 Benoist and Quint \cite{BQ} proved  that if $\mu$ has a moment of order 2, then
\beq \label{var}
\lim_{n \rightarrow \infty} \frac 1 n  \E (S_n^2) = s^2>0 \, ,
 \eeq
  and, for any $t \in {\mathbb R}$, 
$$
\lim_{n \rightarrow \infty}  \sup_{\|x\|=\|y\|=1}  \left|\p \left( \log |\langle A_n x, y \rangle |-n\lambda_\mu\le 
t\sqrt n\right)-\phi(t/s)\right |=0 \, ,
$$
where $\phi$ denotes  the cumulative distribution function of the standard normal distribution. 
\medskip

Given a matrix $g\in GL_d(\R)$ denote by $\lambda_1(g)$ its spectral radius (the greatest modulus of it eigenvalues). Aoun \cite{Aoun} 
proved that if $\mu$ has a moment of order 2, then, for any $t \in {\mathbb R}$, 
$$
\lim_{n \rightarrow \infty}    \left|\p \left( \log (\lambda_1(A_n))-n\lambda_\mu\le 
t\sqrt n\right)-\phi(t/s)\right |=0 \, .
$$

In this paper we provide rates of convergence  in these Central Limit Theorems, if $\mu$ has either an exponential moment, or a polynomial moment of order $p\geq 3$. 

\smallskip

Before giving our main results, les us recall the known results on this subject.
Let $\mu$ be a proximal and strongly irreducible probability measure on ${\mathcal B} (G)$.

If $\mu$ has an exponential moment, then a Berry Esseen bound of order $O(1/\sqrt n)$  for the quantity $\log \|A_n x\|-n \lambda_\mu$ is proved in \cite{LP}.  The same rate is obtained in \cite{CDMP} under a polynomial moment of order 4; in the same paper, the rate $O(\log n / \sqrt n)$ is proved under a moment of order 3. Recently, the rate $O(1/\sqrt n)$ has been obtained in \cite{DKW0} under a moment of order 3, in the special case $d=2$.

If $\mu$ has an exponential moment, then a Berry Esseen bound of order $O(\log n/\sqrt n)$  for the quantity $\log \|A_n \|-n \lambda_\mu$ is proved in \cite{XGL}. The rate $O(1/\sqrt n)$ is obtained in \cite{CDMP} under a polynomial moment of order 4; in the same paper, the rate $O(\log n / \sqrt n)$ is proved under a moment of order 3.

If $\mu$ has an exponential moment, then a Berry Esseen bound of order $O(1/\sqrt n)$  for the quantity $\log |\langle A_n x, y \rangle |-n\lambda_\mu$ has been obtained very recently  by Dinh et al. \cite{DKW} (see also \cite{XGL1} for a more precise statement). This improves on the rate $O(\log n/\sqrt n)$ of Item 1 of Theorem  \ref{matrix-coefficient-exp} below (note that the preprint \cite{DKW} was available on arxiv after this note was submitted).

If $\mu$ has an exponential moment, then a Berry Esseen bound of order $O(\log n/\sqrt n)$  for the quantity $\log (\lambda_1(A_n))-n \lambda_\mu$ is proved in \cite{XGL}.

As we can see, with regard to the Berry-Esseen type bounds for the four quantities described above, the main questions which remain to be treated concern the case of polynomial moments. In particular, it would be interesting to see if the existing moment conditions are optimal (with regard to the rates obtained), and also to propose bounds in  the case where $\mu$ has a polynomial moment of order between 2 and 3.

\section{The case of matrix coefficients}
\setcounter{equation}{0}

\begin{Theorem} \label{matrix-coefficient-exp}
Let $\mu$ be a proximal and strongly irreducible probability measure on ${\mathcal B} (G)$. 
\begin{enumerate}
\item Assume that $\mu$ has an exponential moment, and let $s>0$ be defined by \eqref{var}.
Then there exists a positive constant $K$ such that, for any integer $n\geq 2$, 
\beq\label{case:expo}
\sup_{\|x\|=\|y\|=1} \sup_{t \in {\mathbb R}} \left|\p \left( \log |\langle A_n x, y \rangle |-n\lambda_\mu\le 
t\sqrt n\right)-\phi(t/s)\right |\le \frac{K \log n}{\sqrt n}\, . 
\eeq
\item Assume that $\mu$ has a polynomial moment of order $p\geq 3$ and let $s>0$ be defined by \eqref{var}. Then there exists a positive constant $K$ such that, for any integer $n\geq 2$,
\beq 
\sup_{\|x\|=\|y\|=1} \sup_{t \in {\mathbb R}} \left|\p \left( \log |\langle A_n x, y \rangle |-n\lambda_\mu\le 
t\sqrt n\right)-\phi(t/s)\right |\le \frac{K }{n^{(p-1)/2p}}\, . 
\eeq
\end{enumerate}
\end{Theorem}

%\begin{Remark}
%Very recently Dinh, Kaufmann and Wu \cite{DKW} obtained the bound \eqref{case:expo}  without the extra logarithmic %factor.
%\end{Remark}

The proof of this theorem is based on Berry-Esseen estimates for $\log \|A_n x\|-n \lambda_\mu$ (given in \cite{LP} and \cite{CDMP}), and on  the following elementary lemma (see lemma 5.1 in \cite{H} for a similar result):

\begin{Lemma}\label{basic-lemma}
Let $(T_n)_{n\in \BBN}$ and $(R_n)_{n\in \BBN}$ be two sequences of random variables. Assume that there exist three sequences of positive numbers $(a_n)_{n\in \BBN}$, $(b_n)_{n\in \BBN}$ and $(c_n)_{n\in \BBN}$ going to infinity as $n \rightarrow \infty$,  and a positive constant $s$ such that, for any integer $n$,
\begin{equation*}
\sup_{t\in \BBR}\left|\P(T_n\le  t\sqrt n)-\phi(t/s)\right| \le \frac{1}{a_n}\, ,
\quad \text{and} \quad 
\p(|R_n| \ge \sqrt {2\pi n} s/b_n)\le \frac{1}{c_n}\, .
\end{equation*}
Then, for any integer $n$,
\begin{equation*}
\sup_{t\in \BBR}\left|\P(T_n+R_n\le  t\sqrt n)-\phi(t/s)\right|\le \frac{1}{a_n} +  \frac{1}{b_n} +  \frac{1}{c_n}\, .
\end{equation*}
\end{Lemma}

\noindent {\bf Proof of Lemma \ref{basic-lemma}.} 
Recall that $\phi$ is $1/\sqrt{2\pi}$-Lipschitz. We have 
\begin{align*}
\p(T_n+R_n \le t\sqrt n)  &  \le \p\left(T_n -\sqrt {2\pi n} s/b_n \le t\sqrt n,\, -R_n\le  \sqrt {2\pi n} s/b_n\right)+
\p\left(-R_n \ge  \sqrt {2\pi n} s/b_n\right)\\  &  \le  \p\left(T_n -\sqrt {2\pi n} s/b_n \le t\sqrt n\right) +\p\left(-R_n \ge  
\sqrt{2\pi  n} s/b_n\right)\, .
\end{align*}
Hence
\begin{align*}
\p(T_n+R_n \le t\sqrt n) -\phi(t/s)
&\le \frac{1}{a_n} +|\phi (t/s+\sqrt{2\pi }/b_n)-\phi(t/s)|+ \frac{1}{c_n} \\
& \le \frac{1}{a_n} +  \frac{1}{b_n} +  \frac{1}{c_n}\,  . 
\end{align*}
The lower bound may be proved similarly, by noting that 
\begin{align*}
\p(T_n+ \sqrt {2\pi n} s/b_n \le t\sqrt n)  - \p(R_n\geq \sqrt {2\pi n} s/b_n) &\leq \p(T_n+R_n \le t\sqrt n, \, R_n\leq \sqrt {2\pi n} s/b_n)\\ 
&\leq \p(T_n+R_n \le t\sqrt n) \, . \quad \quad \quad \quad \quad \quad \quad  \quad \quad \quad \square
\end{align*}

\medskip

\noindent {\bf Proof of Item 1 of Theorem \ref{matrix-coefficient-exp}.} The proof follows the steps used in Section 8.3 of \cite{CDJ}.  We shall need some notations. 
For every $\bar x,\, \bar y\in X$, let
$$
d(\bar x, \bar y) := \frac{\|x \wedge y\|}{\|x\| \|y\|}\, ,
$$
where $\wedge$ stands for the exterior product, see e.g. \cite[page 61]{BL},  for the definition and some
properties. Then, $d$ is a metric on $X$.
Let also
\begin{equation}\label{delta}
\delta(\bar x,\bar y):= \frac{|\langle x, y\rangle|}
{\|x\|\, \|y\|}\, .
\end{equation}
Recall that the function $\delta$ is linked to the distance $d$ on $X$ by the following: For every $\bar x, \, \bar y\in X$,

\begin{equation}\label{delta-distance}
\delta^2(\bar x, \bar y)= 1-d^2(\bar x, \bar y)\, .
\end{equation}
\medskip

We shall also need the following result due to Guivarc'h \cite{G} (see Theorem 14.1
in \cite{BQ-Book}):

\begin{Proposition}\label{Guivarch}
Let $\mu$ be a proximal and strongly irreducible probability measure on ${\mathcal B} (G)$. Assume that $\mu$ has an exponential moment. Then, 
there exists $\eta>0$, such that
$$
\sup_{\bar y\in X}\int_X \frac1{\delta(\bar x, \bar y)^\eta }\, d\nu(\bar x)\, <\infty\, . 
$$ 
\end{Proposition}

\medskip

We start with the identity, for $\|x\|=\|y\|=1$,
\begin{align*}
\log |\langle A_n x, y \rangle | &=\log \|A_n x\|  + 
\log \frac{|\langle A_n x, y \rangle |}{\|A_n x\| \|y\|}\\
&= \log \|A_n x\|  + \log \delta(A_n \cdot \bar x, \bar y)\, .
\end{align*}
We shall then apply Lemma \ref{basic-lemma} to $T_n=\log \|A_n x\|-n \lambda_\mu$ and $R_n=\log \delta(A_n \cdot \bar x, \bar y)$. Since $\mu$ has an exponential moment, we know from \cite{LP} that we can take $a_n= C \sqrt n$ in Lemma \ref{basic-lemma}.

In view of Lemma \ref{basic-lemma}, we see  that Theorem 
\ref{matrix-coefficient-exp} will be proved if we can show that there exist
$\tau, K>0$ such that (recall that $\delta (\cdot, \cdot) \le 1$)
\begin{equation}\label{bncn}
\p\left (\left |\log \delta ( A_n \cdot \bar x,\bar y)\right| >\tau \log n\right) 
=\p( \delta ( A_n \cdot \bar x,\bar y)< n^{-\tau})
\le \frac{K}{\sqrt n}\, , 
\end{equation}
which means that  the sequences  $(b_n)_{n\in \BBN}$ and $(c_n)_{n\in \BBN}$ are such that $b_n=\sqrt{2\pi n} s /(\tau \log n)$ and $c_n= \sqrt n /K$.

Recall the identity  \eqref{delta-distance}. As in \cite{CDJ}, we have, using that $d(\cdot, \cdot)\le 1$,
\begin{align} \label{B1}
\delta^2( A_n \cdot \bar x,\bar y)  &   = 1- d^2( A_n \cdot \bar x,\bar y)
\ge 1-\big(d( A_n \cdot \bar x,W_n)+d(W_n, \bar y)\big)^2 \nonumber\\ 
  &  \geq \delta^2(W_n, \bar y)-d^2( A_n \cdot \bar x,W_n )-2d( A_n \cdot \bar x,W_n )d(W_n, \bar y)  \nonumber\\
  &  \ge \delta^2(W_n, \bar y)-3d( A_n \cdot \bar x,W_n )\, .
\end{align}

Hence, to prove  \eqref{bncn}, it suffices to prove that there exist $\tau,K>0$ such that,
\begin{equation}\label{last}
\p\left(\delta^2(W_n, \bar y)< n^{-2\tau} + 3d( A_n \cdot \bar x,W_n )\right)
\le \frac{K }{\sqrt n}\, .
\end{equation}
Now, since $\mu$ has a polynomial moment of order 3, by Lemma 6 of \cite{CDJ}, there exists $\ell >0$, such that 
$$
\p\left(d( A_n \cdot \bar x,W_n ) \ge {\rm e}^{-\ell n}\right)\le 
\frac{C}{n}\, 
$$
(in fact this estimate remains true as soon as $\mu$ has a polynomial moment of order 2, via a monotonicity argument).

Hence, for $n$ large enough (such that 
$3{\rm e}^{-\ell n}\le n^{-2\tau}$), we have 
$$
\p\left(\delta^2(W_n, \bar y)< n^{-2\tau} + 3d( A_n \cdot \bar x,W_n )
\right)
\le \p\left(\delta^2(W_n, \bar y)< 2n^{-2\tau}\right) + \frac{C}n\, .
$$

On another hand, by  Markov's inequality, 
since $W_n$ has law $\nu$,
$$
\p\left(\delta^2(W_n, \bar y )< 2n^{-2\tau}\right)= 
\nu \left\{\bar x\in X\, :\, \frac1{\delta^2(\bar x, \bar y)}>  \frac {n^{2\tau} } 2\right\}
\le  \frac{2^{\eta / 2} } { n^{\eta \tau}} \sup_{\bar y\in X} \int_X \frac1{\delta(\bar x,\bar y)^\eta }\, d\nu(\bar x)\, ,
$$
and \eqref{last} follows  from Proposition \ref{Guivarch} by taking $\tau=\frac1{2\eta}$.  \hfill $\square$

\medskip

\noindent {\bf Proof of Item 2 of Theorem \ref{matrix-coefficient-exp}.} The proof follows the lines of that of Item 1. Instead of Proposition \ref{Guivarch}, we shall use  the following result due to Benoist and Quint (see Propositon 4.5 in \cite{BQ}):

\begin{Proposition}\label{who}
Let $\mu$ be a proximal and strongly irreducible probability measure on ${\mathcal B} (G)$. Assume that $\mu$ has a polynomial moment of order $p>1$. Then
$$
\sup_{\bar y\in X}\int_X  |\log \delta(\bar x, \bar y) |^{p-1}\, d\nu(\bar x)\, <\infty\, . 
$$ 
\end{Proposition}

We shall then apply Lemma \ref{basic-lemma} to $T_n=\log \|A_n x\| - n \lambda_\mu$ and $R_n=\log \delta(A_n \cdot \bar x, \bar y)$. Since $\mu$ has a moment of order 3, we know from \cite{CDMP} that we can take $a_n= C \sqrt n/\log n$   in Lemma \ref{basic-lemma} (and even $a_n= C \sqrt n$ if $p\geq 4$).

In view of Lemma \ref{basic-lemma}, we see that Theorem 
\ref{matrix-coefficient-exp} will be proved if we can show that there exists 
$K>0$ such that 
\begin{equation}\label{bncn2}
\p\left (\left|\log \delta ( A_n \cdot \bar x,\bar y)\right| > n^{1/2p}\right) 
\le \frac{K}{ n^{(p-1)/2p}}\, , 
\end{equation}
which means that  the sequences  $(b_n)_{n\in \BBN}$ and $(c_n)_{n\in \BBN}$ are such that $b_n=\sqrt{2 \pi n} s /n^{1/2p}$ and $c_n= n^{(p-1)/2p}/K$.

Starting again from  \eqref{B1}, we see that it suffices to prove that there exists $K>0$ such that,
\begin{equation}\label{last2}
\p\left(\delta^2(W_n, \bar y)<{\rm e}^{-2n^{1/2p}}+ 3d( A_n \cdot \bar x,W_n )\right)
\le  \frac{K}{ n^{(p-1)/2p}}\, .
\end{equation}

Proceeding as in the proof of Theorem \ref{matrix-coefficient-exp}, we deduce that, for  $n$ large enough (such that 
${\rm e}^{-\ell n}\le {\rm e}^{-2n^{1/2p}}$), we have 
$$
\p\left(\delta^2(W_n, \bar y)< {\rm e}^{-2n^{1/2p}} + 3d( A_n \cdot \bar x,W_n )
\right)
\le \p\left(\delta^2(W_n, \bar y)< 4{\rm e}^{-2n^{1/2p}}\right) + \frac{C}n\, .
$$

On another hand, by  Markov's inequality, 
since $W_n$ has law $\nu$, and for $n$ large enough,
\begin{multline*}
\p\left(\delta^2(W_n, \bar y )< 4{\rm e}^{-2n^{1/2p}}\right)= 
\p\left(|\log \delta (W_n, \bar y)|>  n^{1/2p} - \log 2 \right)  \\
= \nu \left\{\bar x\in X\, :\, |\log \delta (\bar x, \bar y)|>  n^{1/2p} - \log 2\right\} \\
\le  \frac{1} {  (n^{1/2p} - \log 2)^{p-1}} \, \sup_{\bar y\in X}\int_X  |\log \delta(\bar x, \bar y) |^{p-1}\, d\nu(\bar x)\, ,
\end{multline*}
and \eqref{last2} follows  from Proposition \ref{who}.  \hfill 
$\square$

\section{The case of the spectral radius}
\setcounter{equation}{0}

We now prove similar results for the spectral radius. Given a matrix $g\in GL_d(\R)$ denote by $\lambda_1(g)$ its spectral radius (the greatest modulus of its  eigenvalues).  The first result (Item 1 of Theorem \ref{spectral-radius-exp} below), assuming an exponential moment for $\mu$,  has been recently proved by Xiao et al. \cite{XGL} (in fact, a stronger result is proved in \cite{XGL}). We state it  only for  completeness.

\begin{Theorem} \label{spectral-radius-exp}
Let $\mu$ be a proximal and strongly irreducible probability measure on ${\mathcal B} (G)$. 
\begin{enumerate}
\item Assume that $\mu$ has an exponential moment, and let  $s>0$ be defined by \eqref{var}.  Then there exists a positive constant $K$ such that, for any integer $n\geq 2$,
\beq 
\sup_{t \in {\mathbb R}} \left|\p \left( \log \lambda_1(A_n)-n\lambda_\mu\le 
t\sqrt n\right)-\phi(t/s)\right |\le \frac{K \log n}{\sqrt n}\, . 
\eeq
\item 
Assume that $\mu$ has a polynomial moment of order $p\geq 3$, and let  $s>0$ be defined by \eqref{var}.  Then there exists a positive constant $K$ such that, for any integer $n\geq 2$,
\beq 
\sup_{t \in {\mathbb R}} \left|\p \left( \log \lambda_1(A_n)-n\lambda_\mu\le 
t\sqrt n\right)-\phi(t/s)\right |\le \frac{K }{n^{(p-1)/2p}}\, . 
\eeq
\end{enumerate}
\end{Theorem}

The proof of Item 2 is based on  a Berry-Esseen estimate for $\log \|A_n \|-n \lambda_\mu$ given in \cite{CDMP}, and on  Lemma \ref{basic-lemma}.

A key ingredient in the proof of Item 1 by Xiao et al. \cite{XGL} is Lemma 14.13 
of \cite{BQ-Book}. 

To prove Item 2, we shall need  a suitable version of 
Lemma 14.13 of \cite{BQ-Book}. The proof of Lemma 14.13 relies on Lemma 14.2 
of \cite{BQ-Book} (of geometrical nature) and on large deviations, yielding to Lemma 14.3. 

We shall need the following consequence of large deviation estimates of Benoist and Quint \cite{BQ} (see also \cite{CDM} for a related results under proximality). 

\begin{Lemma}\label{lemme-BQ}
Let $\mu$ be a  strongly irreducible probability measure on ${\mathcal B} (G)$. Assume that $\mu$ has a polynomial moment of order $p>1$. 
Let $\varepsilon>0$. There exists $C>0$ such that for every $n\in \N$
\begin{gather}
\label{BQ1}\sup_{\|x\|=1} \P \left(\max_{1 \leq k \leq n}|\log \|A_k x\|-k \lambda_\mu|>\varepsilon n \right)
\le \frac{C}{n^{ p-1}}\, ,\\
\label{BQ2} \P\left(\max_{1 \leq k \leq n}|\log \|A_k\|-n \lambda_\mu|>\varepsilon n \right)
\le \frac{C}{n^{ p-1}}\, ;\\
\label{BQ3} \P\left (\max_{1 \leq k \leq n}|\log \|\Lambda^2(A_k)\|-k (\lambda_\mu 
 +\gamma_\mu)|>\varepsilon n \right)
\le \frac{C}{n^{ p-1}}\, .
\end{gather}
\end{Lemma}
\noindent {\bf Remark.} Let us recall that, for any $A \in GL_d({\mathbb R})$, $\Lambda^2(A)$ is the matrix on $\Lambda^2({\mathbb R }^d)$ defined by $\Lambda^2(A)(x \wedge y) = Ax \wedge A y$.  In addition, in  \eqref{BQ3} $,\gamma_\mu$ is the second Lyapunov exponent of $\mu$. With the notations of 
\cite[Section 14]{BQ-Book}, $\lambda_\mu$ is denoted either $\lambda_{1,\mu}$ or $\lambda_1$, while $\gamma_\mu$ is either denoted $\lambda_{2,\mu}$ or $\lambda_2$. 

\medskip

\noindent {\bf Proof of Lemma \ref{lemme-BQ}.}  Let $u_n$ be any of the left-hand side in \eqref{BQ1}, \eqref{BQ2} or \eqref{BQ3}. It follows from Proposition 4.1 and Corollary 4.2 of \cite{BQ} that 
\beq \label{BQmax}
\sum_{n\ge 1}  n^{p-2} u_n<\infty \, .
\eeq
In fact,  in \cite{BQ}, \eqref{BQmax} is proved for  $u_n$ defined without the maximum over $k \in \{1, \ldots , n\}$ under the probability. However, it is easy to see that the maximum over $k$ can be added: it suffices to follow the proof of Theorem 2.2 of 
\cite{BQ} with obvious changes, and to use a maximal version of Burkholder's inequality for martingales. Now,  once \eqref{BQmax} has been proven, it is easy to infer (via a monotonicity argument) that  \eqref{BQ1}, \eqref{BQ2} and \eqref{BQ3} are satisfied. 
\hfill $\square$

\medskip

Using Lemma \ref{lemme-BQ} one can reproduce the proof of Proposition 14.3 of \cite{BQ-Book} to prove the following version of it. 

\begin{Lemma}\label{lemme-aux}
Let $\mu$ be a  strongly irreducible and proximal probability measure on ${\mathcal B} (G)$. Assume that $\mu$ has a polynomial moment of order $p>1$. Then, the estimates (14.5), (14.6), (14.7) and (14.8) 
of \cite{BQ-Book} hold with $1- \frac{C}{n^{p-1}}$ in the right-hand side instead of 
$1-{\rm e}^{-cn}$.
\end{Lemma}

\medskip

Lemma \ref{lemme-aux} implies the next result. 

\begin{Lemma}\label{control}
Let $\mu$ be a  strongly irreducible probability measure on ${\mathcal B} (G)$. Assume that $\mu$ has a polynomial moment of order $p>1$. For every $\varepsilon>0$, there exist $C>0$ and $\ell_0>0$ such that for every $\ell_0\le \ell 
\le n$, 
$$
\P(\log(\lambda_1(A_n))-\log \|A_n\|\ge -\varepsilon \ell)\ge 1-\frac{C}{\ell^{p-1}}\, .
$$
\end{Lemma}

\noindent {\bf Proof of Lemma \ref{control}.} The lemma is a version of Lemma 14.13 of \cite{BQ-Book} with the following difference:  Lemma 14.13 holds under an exponential moment while in Lemma \ref{control} we assume polynomial moments. 

\smallskip

Now, it happens that there is a small gap in the proof of Lemma 14.13 of \cite{BQ-Book} which can be  easily fixed thanks to a slight modification of the original argument.

\medskip

One of the steps in the proof of Lemma 14.13 consists in 
proving  that the property (14.38) is true on an exponentially small set (see the end of page 233 of \cite{BQ-Book} for the definition of an exponentially small set).  A second  step of the proof  consists in proving the equivalence of the fact that 
the properties  (14.38) and (14.43)  of \cite{BQ-Book} are true on an exponentially small set . 

\medskip

The problem then comes from the fact that it does not seem possible to deduce straightforwardly from (14.7) that the property (14.43) is true on an exponentially small set, as mentioned in 
\cite{BQ-Book}.  Yet the weaker property \eqref{ine-BQ} below follows from (14.7). Notice that since we prove below that the property (14.38) is true on an exponentially small set, from the above mentioned equivalence, it will follow that the property (14.43) is also true on an exponentially small set.

\medskip

We choose to explain how to fix the proof of the original Lemma 14.13. Then, the proof of our Lemma \ref{control} may be done similarly, using our Lemma \ref{lemme-aux} instead of Lemma 14.3 of \cite{BQ-Book}.

\smallskip

From (14.7) of \cite{BQ-Book} it follows that, with the notations of \cite{BQ-Book}, for every $n\ge n_0$
\begin{equation}\label{ine-BQ}
\mu^{\otimes n}\Big( \Big \{(b_1, \ldots , b_n)\in G^n \, :\, 
\delta(x^M_{b_n\cdots b_{[n/2]+1}},y^m_{b_{[n/2]}\cdots b_1})\ge {\rm e}^{-\varepsilon [n/2]} \Big \}\Big)\ge 1-{\rm e}^{-c [n/2]}\, .
\end{equation}
Using (14.39), (14.40), (14.41) and (14.42), this yields 
that 
\begin{equation}\label{ine-BQ-2}
\mu^{\otimes n}\Big( \Big \{(b_1, \ldots , b_n)\in G^n \, :\, 
\delta(x^M_{b_n\cdots b_{1}},y^m_{b_{n}\cdots b_1})\ge {\rm e}^{-\varepsilon \ell} \Big \}\Big)\ge 1-{\rm e}^{-c \ell}
\qquad \forall [n/2]\le \ell \le n\,, 
\end{equation}
for some $c>0$ that may differ from the above one (and from the other $c$'s below).

\medskip

Let $\ell_0\le \ell <[n/2]$, with $\ell_0\ge n_0$, where $n_0$ is such that Lemma 14.3 be true. 

\smallskip

By (14.6) of \cite{BQ-Book}, we have 
\begin{equation}
\label{ineq1}\mu^{\otimes n}\Big(\Big \{(b_1, \ldots , b_n)\in G^n \, :\, 
d(x^M_{b_n\cdots b_{n-\ell}},b_n\cdots b_1x_0)\le {\rm e}^{-(\lambda_{1,\mu}-\lambda_{2,\mu}-\varepsilon)\ell}\Big\}\Big)\ge 1-{\rm e}^{-c \ell}\, . 
\end{equation}

By (14.7), we have
\begin{gather}
\label{ineq2}\mu^{\otimes n}\Big( \Big \{(b_1, \ldots , b_n)\in G^n \, :\, 
\delta (x^M_{b_n\cdots b_{n-\ell}},y^m_{b_{[n/2]}\cdots b_1})\ge {\rm e}^{-\varepsilon \ell}\Big\}\Big)\ge 1-{\rm e}^{-c \ell}\, ,
\end{gather}
where we used that $b_n\cdots b_{n-\ell}$ and $b_{[n/2]}
\cdots b_1$ are independent since $n-\ell >[n/2]$.

\medskip

Using the fact that (14.39), (14.41) and (14.42) are true except on an exponentially small set, combined with 
\eqref{ineq1} and \eqref{ineq2}, we infer that 
\begin{equation}\label{ine-BQ-3}
\mu^{\otimes n}\Big( \Big \{(b_1, \ldots , b_n)\in G^n \, :\, 
\delta(x^M_{b_n\cdots b_{1}},y^m_{b_{n}\cdots b_1})\ge {\rm e}^{-\varepsilon \ell} \Big \}\Big)\ge 1-{\rm e}^{-c \ell}
\qquad \forall 1\le \ell <[n/2]\,.
\end{equation}
Combining \eqref{ine-BQ-2} and \eqref{ine-BQ-3}, 
 we see that the property (14.38) is 
true on an exponentially small set. 

\medskip
Then, the proof of Lemma 14.13 may be finished as in 
\cite{BQ-Book}, combining Lemma 14.14 with the facts that 
the properties (14.37) and (14.38) are true on  exponentially small sets. \hfill $\square$

\medskip

\noindent {\bf Proof of Item 2 of Theorem \ref{spectral-radius-exp}.} 
We shall  apply Lemma \ref{basic-lemma} to $T_n=\log \|A_n \|-n \lambda_\mu$ and $R_n=\log(\lambda_1(A_n))-\log \|A_n\|$. Since $\mu$ has a moment of order 3, we know from \cite{CDMP} that we can take $a_n= C \sqrt n/\log n$  in Lemma \ref{basic-lemma} (and even $a_n= C \sqrt n$ if $p\geq 4$).

In view of Lemma \ref{basic-lemma}, we see that Item 2 of Theorem 
\ref{spectral-radius-exp} will be proved if we can show that there exists 
$K>0$ such that 
\begin{equation}\label{bncn4}
\p\left (\left|\log(\lambda_1(A_n))-\log \|A_n\|\right| > n^{1/2p}\right) 
\le \frac{K}{ n^{(p-1)/2p}}\, , 
\end{equation}
which means that  the sequences  $(b_n)_{n\in \BBN}$ and $(c_n)_{n\in \BBN}$ are such that $b_n=\sqrt{2 \pi n} s /n^{1/2p}$ and $c_n= n^{(p-1)/2p}/K$.

Recall that 
 $\lambda_1(g)\le \|g\|$ for every $g\in GL_d(\R)$. Hence  \eqref{bncn4} follows from Lemma \ref{control} by taking $\varepsilon=1$ and  $\ell=n^{1/2p}$. \hfill $\square$

\end{document}